\documentclass[12pt]{article}
\usepackage{amsfonts}
\usepackage{amsmath,amsthm,amscd,amssymb,
indentfirst,bm,makeidx}
\usepackage{mathrsfs}

\usepackage[dvips]{graphicx}
\usepackage[all]{xy}
\usepackage{graphicx}
\usepackage{fancyhdr}
\usepackage{CJK}
\usepackage{bbm}

\def\y{\begin{eqnarray*}}
\def\bd{\begin{description}}
\def\ey{\end{eqnarray*}}
\def\ebd{\end{description}}

\begin{document}

\setlength{\baselineskip}{16pt}

\title{Asymptotic behavior of time periodic solutions for extended Fisher-Kolmogorov equations with delays \thanks{Research supported by NNSF of China (12061063, 12061062, 11771428), Project of NWNU-LKQN2019-3 and China Scholarship Council (201908625016).}
}
\author{Pengyu Chen $^\textrm{a,}$\thanks{Corresponding author. E-mail addresses: chpengyu123@163.com (P. Chen), lanyu9986@126.com (X. Zhang), zhzt@math.ac.cn (Z. Zhang).},
Xuping Zhang $^\textrm{a}$, Zhitao Zhang $^\textrm{b,c}$\\
{\small $^\textrm{a}$ Department of Mathematics, Northwest Normal University, Lanzhou 730070, P.R. China}\\
{\small $^\textrm{b}$  Academy of Mathematics and Systems Science, Chinese Academy}\\
{\small of Sciences, Beijing 100190, P. R. China   }\\
{\small $^\textrm{c}$  School of Mathematical Sciences, University of Chinese Academy}\\
{\small  of Sciences, Beijing 100049, P. R. China }}

\date{}

\maketitle
\begin{abstract}
\setlength{\baselineskip}{14pt} In this paper, we investigate the global existence, uniqueness and asymptotic stability of time periodic classical solution for a class of extended Fisher-Kolmogorov equations with delays and  general nonlinear term. We establish a general framework to investigate the asymptotic behavior of time periodic solutions for nonlinear extended Fisher-Kolmogorov equations with  delays and general nonlinear function, which will provide an effective way to deal with such kinds of problems. The discussion is based on the theory of compact and analytic operator semigroups and maximal  regularization method.
\vspace{8pt}

 \vspace{3pt}

\end{abstract}

\setlength{\baselineskip}{16pt}
\vskip10mm
\section{Introduction and main results}

The extended Fisher-Kolmogorov (EFK) equation
$$
\frac{\partial u}{\partial t}+\gamma \frac{\partial^4 u}{\partial x^4}-\frac{\partial^2 u}{\partial x^2}-u+u^3=0,\quad \gamma>0 \eqno(1.1)
$$
was proposed in 1987 by Coullet, Elphick, and Repaux \cite{cer87} and in 1988
by Dee and van Saarloos \cite{ds88} as a generalization of the classical Fisher-Kolmogorov
equation which was firstly propounded by Fisher and Kolmogorov in 1937. The extended Fisher-Kolmogorov equation are a family of models arising in population dynamics problems, cancer modelling, chemical kinetics, the description of propagating crystallization/polymerization fronts, geochemistry and
many other fields. These equations
do not admit a Lagrangian density depending on the field $u$ and thus the
variational formulation for the effective particle parameters cannot be written
in the usual way. Therefore, substantial attention has been focused on the steady-state equation
$$
-\gamma u''''+u''+u-u^3=0,\quad \gamma>0 \eqno(1.2)
$$
corresponding to EFK equation (1.1), see \cite{de14,gst05,pt97,sv02} and references therein for more comments and citations.\vskip3mm

Recently, Danumjaya and Paniuse employing the Galerkin finite element approximation method and the orthogonal cubic spline collocation method studied the existence, uniqueness and regularity of EFK equation (1.1) in \cite{dp05} and \cite{dp06}. In 2011, by using a Crank-Nicolson type
finite difference scheme and the method of Lyapunov functional, Khiari and Omrani \cite{ko11}
studied the existence of approximate solutions for the following extended Fisher-Kolmogorov equation in two space dimension with Dirichlet boundary conditions
$$\left\{\begin{array}{ll}
u_t+\gamma \Delta^2u-\Delta u-u+u^3=0,\quad \textrm{in} \; (0,T]\times \Omega,\\[6pt]
u=0,\quad \Delta u=0,\quad (t,x,y)\in (0,T]\times \partial\Omega,\\[6pt]
u(0,x,y)=u_0(x,y),\quad \textrm{in} \;  \Omega,
 \end{array} \right.\eqno(1.3)
$$
where $\Omega$ is a bounded domain in $\mathbb{R}^2$ with boundary $\partial\Omega$, $\gamma>0$ is a constant.\vskip3mm

On the other hand, evolution equations with delays
have attracted increasing attention in past twenty years and the existence or attractivity of periodic solutions for evolution equations with delays have been considered by several authors, see \cite{bz91,b85,ccn12,h88,h81,l11,ll06,o11,w96,xa92,zlw16} and references listed therein for more comments and citations. Most of these results are established by applying
semigroup theory \cite{ccn12,h88,h81,l11,o11,w96,xa92}, corresponding fixed point theorems \cite{ccn12,l11,ll06,zlw16}, coincidence degree theory \cite{h88} and so on. Recently, Liu and Li \cite{ll06} obtained the existence
of periodic solutions for a class of  parabolic evolution equations with delay by utilizing a Schaefer type theorem, which extend the corresponding results of Burton and  Zhang \cite{bz91}. Latter, in 2008,  by using the method of constructing some suitable Lypunov functionals and establishing an a priori bound for all possible periodic solutions, Zhu, Liu and Li \cite{zll08} investigated the existence, uniqueness and global attractivity of time periodic solutions for the following
one-dimensional parabolic evolution equation with delays

$$
\left\{\begin{array}{ll}
\frac{\partial}{\partial t}u(t,x)=\frac{\partial^2 }{\partial x^2}u(t,x)+au(t,x)+g(t,x)\\[6pt]
\qquad\qquad\quad+f(u(t-\tau_1,x),\cdots,u(t-\tau_n,x)),\quad (t,x)\in \mathbb{R}\times (0,1),\\[12pt]
u(t,0)=u(t,1)=0,\quad t\in \mathbb{R},
 \end{array} \right.\eqno(1.4)
$$
which is usually used to model some process of biology, where $a\in \mathbb{R}$, $f: \mathbb{R}^n\rightarrow \mathbb{R}$ is locally Lipschitz continuous, $g: \mathbb{R}\times[0, 1]\rightarrow \mathbb{R}$ is H\"{o}lder continuous
and $g(t,x)$ is $\omega$-periodic in $t$, $\tau_1$, $\tau_2$, $\cdots$, $\tau_n$ are positive constants.
In addition, the dynamical characteristics (including stable,
unstable, attract, oscillatory and chaotic behavior) of differential
equations have become a subject of intense research activities. For
the details of this field, we refer the reader to the monographs of
Burton \cite{b85}, Hale \cite{h88} and the papers of Caicedo, Cuevasa,
Mophoub and N'Gu\'{e}r\'{e}kata \cite{ccn12}, Chen and Guo \cite{cg02}, Li and Wang \cite{lw10} and Wang, Liu and Liu \cite{wll05}. As far as we know, no work
has been done for the asymptotic behavior of time periodic solutions for the extended Fisher-Kolmogorov equations. This is an interesting and important problem that needs to be solved. Also, it is one of motivations of this paper.\vskip3mm

Naturally, the delay also occurs in the model of EFK extended Fisher-Kolmogorov equations. However, to our best knowledge, up until now the time periodic solutions for extended Fisher-Kolmogorov equations with delays have not
been considered in the literature. Motivated by the above consideration, in this paper, we are concerned with the existence, uniqueness and asymptotic behavior of time $\omega$-periodic classical solutions for the following extended Fisher-Kolmogorov (EFK) equations with delays and general nonlinear term of the form
$$
\left\{\begin{array}{ll}
\frac{\partial}{\partial t}u(t,x)+\gamma \frac{\partial^4}{\partial x^4}u(t,x)-\frac{\partial^2 }{\partial x^2}u(t,x)-u(t,x)=g(t,x)\\[6pt]
\qquad\qquad\qquad\qquad\quad+f(u(t-\tau_1,x),\cdots,u(t-\tau_n,x)),\;\textrm{in}\; \mathbb{R}\times (0,1),\\[12pt]
u(t,0)=u(t,1)=0,\quad u_{xx}(t,0)=u_{xx}(t,1)=0,\quad t\in \mathbb{R},
 \end{array} \right.\eqno(1.5)
$$
where $\gamma>0$ is a constant, $f: \mathbb{R}^n\rightarrow \mathbb{R}$ is a nonlinear continuous function, $g: \mathbb{R}\times[0, 1]\rightarrow \mathbb{R}$ is continuous
and $g(t,x)$ is $\omega$-periodic in $t$, $\tau_1$, $\tau_2$, $\cdots$, $\tau_n$ are positive constants.\vskip3mm

In this paper, by defining a positive definite selfadjoint operator $A$, which generates a compact semigroup $T(t)$ $(t\geq0)$  in Hilbert space $H$, we first transfer the extended Fisher-Kolmogorov equations with delays (1.5) into the abstract form for a class of nonlinear evolution equation in the frame of Hilbert space $H$, and then apply corresponding fixed point
theorems, the theory of compact operator semigroups and nonlinear analysis
theory to discuss the existence and uniqueness of $\omega$-periodic mild solutions for abstract nonlinear evolution equation. Further, by applying the maximal regularity of linear evolution equations with
positive definite operator combined with the regularization
method via the theory of analytic semigroups, we proved the existence and uniqueness of time $\omega$-periodic classical solution for extended Fisher-Kolmogorov equations with delays (1.5). In addition, based on the uniqueness of time $\omega$-periodic classical solution,  we obtained the global asymptotic stability of time $\omega$-periodic classical solution for extended Fisher-Kolmogorov equations with delays (1.5) by using the exponentially stability of analytic semigroup $T(t)$ $(t\geq0)$ and an integral inequality of Bellman type with delays.\vskip3mm
The main results of this paper are as follows:
\vskip3mm\noindent\textbf{Theorem 1.1.} \emph{Assume that $f: \mathbb{R}^n\rightarrow \mathbb{R}$ is locally Lipschitz continuous, $g: \mathbb{R}\times[0, 1]\rightarrow \mathbb{R}$ is H\"{o}lder continuous
and $g(t,x)$ is $\omega$-periodic in $t$. If the following conditions
\begin{itemize}
\item [(H1)] There exist positive constants $\beta_1$, $\beta_2$, $\cdots$, $\beta_n$ and $K$ such that
$$|f(\xi_1,\cdots,\xi_n)+g(t,x)|\leq \sum_{k=1}^{n}\beta_k|\xi_k|+K\quad \textrm{for}\; x\in [0,1],\;(\xi_1,\cdots,\xi_n)\in
\mathbb{R}^n;$$
\item [(H2)] $\sum\limits_{k=1}^{n}\beta_k<\gamma\pi^4+\pi^2-1$,
\end{itemize}
hold, then EFK equation (1.5) has at least one
time $\omega$-periodic classical solution $u\in C^{1,2}(\mathbb{R}\times[0,1])$.}\vskip3mm

If we  strengthen condition (H1), then we have the following uniqueness result of time $\omega$-periodic classical solution for EFK equation (1.5).
\vskip3mm\noindent\textbf{Theorem 1.2.} \emph{Assume that $f: \mathbb{R}^n\rightarrow \mathbb{R}$ is locally Lipschitz continuous, $g: \mathbb{R}\times[0, 1]\rightarrow \mathbb{R}$ is H\"{o}lder continuous
and $g(t,x)$ is $\omega$-periodic in $t$. If the following condition
\begin{itemize}
\item [(H3)] There exist positive constants $\beta_1$, $\beta_2$, $\cdots$, $\beta_n$ such that
$$|f(\xi_1,\cdots,\xi_n)-f(\eta_1,\cdots,\eta_n)|\leq \sum_{k=1}^{n}\beta_k|\xi_k-\eta_k|\quad \textrm{for}\; (\xi_1,\cdots,\xi_n),\,(\eta_1,\cdots,\eta_n)\in
\mathbb{R}^n,$$
\end{itemize}
and condition (H2) hold, then EFK equation (1.5) has a unique
time $\omega$-periodic classical solution $u\in C^{1,2}(\mathbb{R}\times[0,1])$.}\vskip3mm

If we  strengthen condition (H2), then we can obtain the global asymptotic stability of time $\omega$-periodic classical solution for EFK equation (1.5).
\vskip3mm\noindent\textbf{Theorem 1.3.} \emph{Assume that $f: \mathbb{R}^n\rightarrow \mathbb{R}$ is locally Lipschitz continuous and $g:\mathbb{R}_+ \times[0, 1] \rightarrow \mathbb{R}$ is H\"{o}lder continuous. If  condition (H3) and the following condition
\begin{itemize}
\item [(H2)$'$] $\sum\limits_{k=1}^{n}\beta_k e^{(\gamma\pi^4+\pi^2-1)\tau_k}<\gamma\pi^4+\pi^2-1$,
\end{itemize}
hold, then EFK equation (1.5) has a unique
time $\omega$-periodic classical solution $\bar{u}\in C^{1,2}(\mathbb{R}\times[0,1])$ and it is globally asymptotically stable.}\vskip3mm

The rest of this paper is organized as follows: In the following section we first introduce some notations and preliminaries which
are used throughout this paper. Especially, the extended Fisher-Kolmogorov equations with delays (1.5) is transformed into an abstract nonlinear evolution equation in a Hilbert space $H$.  In section
3 we prove the global existence and uniqueness  of
time $\omega$-periodic classical solutions for extended Fisher-Kolmogorov equations with delays (1.5) (Theorems 1.1 and 1.2).
In the last
section, we prove the global asymptotic stability   of time $\omega$-periodic classical solution for extended Fisher-Kolmogorov equations with delays (1.5) (Theorem 1.3).
\vskip20mm

\section{Preliminaries}
\vskip5mm
Let $H=L^2([0,1],\mathbb{R})$ be a real Hilbert space with the $L^2$-norm $\|\cdot\|_2$ defined by
$$\|u\|_2=\Big(\int_0^1|u(x)|^2dx\Big)^{\frac{1}{2}},\qquad \forall\; u\in L^2([0,1],\mathbb{R})$$
and inner product $\langle\cdot,\cdot\rangle$ defined by
$$
\langle u,v\rangle=\int_0^1u(x)v(x)dx,\qquad \forall\; u,v\in L^2([0,1],\mathbb{R}).
$$ We
define an operator $A$ in Hilbert space $H$ by
$$\aligned
&Au=\gamma \frac{\partial^4 u}{\partial x^4}-\frac{\partial^2 u}{\partial x^2}-u,\\[12pt]
&D(A)=\{W^{4,2}[0,1]\mid u(0)=u(1)=u''(0)=u''(1)=0\}.
\endaligned\eqno(2.1)
$$
From (2.1) it is easy to show that $ D(A)$ is densely defined in $H$.\vskip2mm

Let
$u(t)=u(t,\cdot)$,  $f(u(t-\tau_1),\cdots,u(t-\tau_n))=f(u(t-\tau_1,\cdot),\cdots,u(t-\tau_n,\cdot))$, $g(t)=g(t,\cdot)$. Then the extended Fisher-Kolmogorov (EFK) equation with delays (1.5) can be transformed into the abstract form of delay evolution equation
$$
u'(t)+ Au(t)=f(u(t-\tau_1),\cdots,u(t-\tau_n))+g(t),\quad t\in \mathbb{R},
\eqno(2.2)
$$
in the Hilbert space $H=L^2([0,1],\mathbb{R})$.\vskip3mm
\noindent\textbf{Lemma 2.1.} \emph{The operator $A:D(A)\subset H\rightarrow H$ defined by (2.1) is a symmetric operator.
}\vskip2mm
\noindent\textbf{Proof.}\quad For any $u$, $v\in D(A)$, using integration by parts one gets that
\begin{eqnarray*}
\quad\langle Au,v\rangle&=& \gamma \int_0^1\frac{\partial^4 u(x)}{\partial x^4}v(x)dx-\int_0^1\frac{\partial^2 u(x)}{\partial x^2}v(x)dx-\int_0^1u(x)v(x)dx\\[6pt]
&=&\gamma \int_0^1\frac{\partial^2 u(x)}{\partial x^2}\frac{\partial^2 v(x)}{\partial x^2}dx+\int_0^1\frac{\partial u(x)}{\partial x}\frac{\partial v(x)}{\partial x}dx-\int_0^1u(x)v(x)dx\\[6pt]
&=&\gamma \int_0^1 u(x)\frac{\partial^4 v(x)}{\partial x^4}dx-\int_0^1 u(x)\frac{\partial^2 v(x)}{\partial x^2}dx-\int_0^1u(x)v(x)dx\\[6pt]
&=&\langle u,Av\rangle.\qquad\qquad\qquad\qquad\qquad\qquad\qquad\qquad\qquad
\qquad\qquad\qquad\;(2.3)
\end{eqnarray*}
(2.3) means that the operator $A$ is a symmetric operator. This completes the proof of Lemma 2.1.\hfill$\Box$\vskip3mm
\noindent\textbf{Lemma 2.2.} \emph{ The operator $A:D(A)\subset H\rightarrow H$ defined by (2.1) is a positive definite operator.}\vskip2mm
\noindent\textbf{Proof.}\quad For any $u\in D(A)$, by (2.1) and Poincare inequality which can be find in \cite{te97}, we get that
\begin{eqnarray*}
\quad\langle Au,u\rangle&=& \gamma \int_0^1\frac{\partial^4 u(x)}{\partial x^4}u(x)dx-\int_0^1\frac{\partial^2 u(x)}{\partial x^2}u(x)dx-\int_0^1u(x)u(x)dx\\[6pt]
&=&\gamma \int_0^1\frac{\partial^2 u(x)}{\partial x^2}\frac{\partial^2 u(x)}{\partial x^2}dx+\int_0^1\frac{\partial u(x)}{\partial x}\frac{\partial u(x)}{\partial x}dx-\int_0^1u^2(x)dx\\[6pt]
&=&\gamma \Big\|\frac{\partial^2 u}{\partial x^2}\Big\|_2^2+\Big\|\frac{\partial u}{\partial x}\Big\|_2^2-\|u\|_2^2\\[6pt]
&\geq&(\gamma\pi^4+\pi^2-1)\|u\|_2^2\\[6pt]
&\geq&
\|u\|_2^2,
\end{eqnarray*}
and $\langle Au,u\rangle=0$ if and only if $u=0$. Therefore, $A$ is a positive definite operator. This completes the proof of Lemma 2.2. \hfill$\Box$\vskip3mm
\noindent\textbf{Lemma 2.3.}\emph{ $\mathcal{R}(A)=H$.}\vskip2mm
\noindent\textbf{Proof.}\quad We only need to prove that for any $\phi\in H$ there exist $u\in D(A)$ such that $Au=\phi$. This fact is equivalent to solve the following linear boundary value problem of fourth-order ordinary differential equation
$$
\left\{\begin{array}{ll}
\gamma \frac{d^4 u(x)}{d x^4}-\frac{d^2 u(x)}{d x^2}-u(x)=\phi(x),\quad x\in [0,1],\\[14pt]
u(0)=u(1)=u_{xx}(0)=u_{xx}(1)=0,
 \end{array} \right.\eqno(2.4)
$$
namely
$$
Au=\gamma\Big(-\frac{d^2}{d x^2}+\mu_1\Big)\Big(-\frac{d^2}{d x^2}+\mu_2\Big)u(x)=\phi(x), \eqno(2.5)
$$
where
$$
\mu_{1}=\frac{1+ \sqrt{1+4\gamma}}{2\gamma}, \qquad \mu_{2}=\frac{1- \sqrt{1+4\gamma}}{2\gamma}.\eqno(2.6)
$$
Since $\gamma>0$, it is obvious that  $\mu_{1}>0$.  From (2.6) one gets that
$$
 \mu_{2}+\pi^2=\frac{1- \sqrt{1+4\gamma}}{2\gamma}+\pi^2=\frac{1- \sqrt{1+4\gamma}+2\gamma\pi^2}{2\gamma}>\frac{1+4\gamma- \sqrt{1+4\gamma}}{2\gamma}>0,
$$
which means that $-\pi^2<\mu_{2}<0$. By \cite{l03} we know that the solution of the linear boundary value problem (2.4) can be expressed by
$$
u(x)=\gamma\int_0^1\int_0^1G_1(x,y)G_2(y,z)\phi(z)dz dy,\eqno(2.7)$$
where $G_i(x,y)$ ($i=1$, $2$) is the Green's function of the second order linear boundary value problem
$$
\left\{\begin{array}{ll}
-u''(x)+\mu_{i}u(x)=0,\quad x\in [0,1],\\[10pt]
u(0)=u(1)=0,
 \end{array} \right.\eqno(2.8)
$$
and $G_i(x,y)$ can be expressed by
$$
G_1(x,y)=\left\{\begin{array}{ll}
\frac{\sinh\sqrt{|\mu_{1}|}x\cdot \sinh\sqrt{|\mu_{1}|}(1-y)}{\sqrt{|\mu_{1}|}\sinh\sqrt{|\mu_{1}|}},\quad 0\leq x\leq y \leq 1,\\[12pt]
\frac{\sinh\sqrt{|\mu_{1}|}y\cdot \sinh\sqrt{|\mu_{1}|}(1-x)}{\sqrt{|\mu_{1}|}\sinh\sqrt{|\mu_{1}|}},\quad 0\leq y\leq x \leq 1,
 \end{array} \right.
$$
and
$$
G_2(x,y)=\left\{\begin{array}{ll}
\frac{\sin\sqrt{|\mu_{2}|}x\cdot \sin\sqrt{|\mu_{2}|}(1-y)}{\sqrt{|\mu_{2}|}\sin\sqrt{|\mu_{2}|}},\quad 0\leq x\leq y \leq 1,\\[12pt]
\frac{\sin\sqrt{|\mu_{2}|}y\cdot \sin\sqrt{|\mu_{2}|}(1-x)}{\sqrt{|\mu_{2}|}\sin\sqrt{|\mu_{2}|}},\quad 0\leq y\leq x \leq 1.
 \end{array} \right.
$$
(2.5) and (2.7) mean that for any $\phi\in H$ there exist $u\in D(A)$ such that $Au=\phi$. Therefore,  $\mathcal{R}(A)=H$. This completes the proof of Lemma 2.3. \hfill$\Box$\vskip3mm

Therefore, from Lemmas 2.1, 2.2 and 2.3, we know that the operator $A:D(A)\subset H\rightarrow H$ defined by (2.1) is a positive definite selfadjoint operator and the first eigenvalue of the operator $A$ is $\lambda_1=\gamma\pi^4+\pi^2-1$. Furthermore, by Lemma 2.3 and \cite{l03}  one can easily to prove that the operator $A:D(A)\subset H\rightarrow H$ defined by (2.1) has compact resolvent. Hence, it is well known from \cite{h81,pa83} that the operator $A$ defined by (2.1) is a sectorial operator, and therefore $-A$ generates an analytic and compact semigroup $T(t)$ $(t\geq0)$ in $H$, which is exponentially stable and satisfies
$$
\|T(t)\|\leq e^{-(\gamma\pi^4+\pi^2-1)t},\quad t\geq 0.\eqno(2.9)
$$
\vskip2mm

Next, we give some concepts and conclusions on the fractional powers of
$A$. For $\alpha>0$, $A^{-\alpha}$ is defined by
$$
A^{-\alpha}=\frac{1}{\Gamma(\alpha)}\int_0^{\infty}s^{\alpha-1}T(s)ds,\eqno(2.10)$$
where $\Gamma(\cdot)$ is the Gamma function. $A^{-\alpha}\in
\mathcal {B}(H)$ is injective, and $A^\alpha$ can be defined by
$A^\alpha=(A^{-\alpha})^{-1}$ with the domain
$D(A^\alpha)=A^{-\alpha}(H)$, where $\mathcal {B}(H)$ denotes the Banach space of all linear bounded operators
from $H$ to $H$ endowed with the topology defined by operator norm. For $\alpha=0$, let $A^\alpha=I$. We
endow $D(A^\alpha)$ with an inner product
$\langle\cdot,\cdot\rangle_\alpha=\langle A^\alpha\cdot,A^\alpha\cdot\rangle$. Since $A^\alpha$ is a closed linear operator, it
follows that $(D(A^\alpha),\langle\cdot,\cdot\rangle_\alpha)$ is a Hilbert
space. We denote by $H_\alpha$ the Hilbert space
$(D(A^\alpha),\langle\cdot,\cdot\rangle_\alpha)$. Especially, $H_0=H$ and
$H_1=D(A)$. For $0\leq\alpha<\beta$, $H_\beta$ is densely embedded
into $H_\alpha$ and the embedding $H_\beta\hookrightarrow H_\alpha$
is compact. For the details, we refer to \cite{h81} and \cite{xa92}.\vskip3mm

From \cite[Chapter 4, Corollary 2.5]{pa83}, we know that for any $u_0\in
D(A)$, if the linear function $\varphi$ is continuously differentiable on $\mathbb{R}_+$, then the initial value problem for the linear
evolution equation (LIVP)
$$
\left\{\begin{array}{ll}
  u'(t)+Au(t)= \varphi(t),\quad  t\in \mathbb{R}_+, \\[12pt]
   u(0)=u_0
 \end{array} \right.\eqno(2.11)
$$
has a unique classical solution $u\in C^1((0,+\infty),H)\cap C((0,+\infty),H_1)\cap C(\mathbb{R}_+$, $H)$
expressed by
$$
u(t)=T(t)u_0+\int_0^{t}T(t-s)\varphi(s)ds.\eqno(2.12)
$$
If $u_0\in H$ and $\varphi\in L^1(\mathbb{R}_+,H)$, the function $u$ given by (2.12)
belongs to $C(\mathbb{R}_+,H)$, which is known as a mild solution of the LIVP
(2.11). If a mild solution $u$ of the LIVP (2.11) belongs to
$W^{1,1}(\mathbb{R}_+,H)\cap L^1(\mathbb{R}_+,H_1)$ and satisfies the equation for a.e.
$t\in \mathbb{R}_+$, we call it a strong solution. By \cite[Chapter 4, Corollary 2.10]{pa83}, we know that for any $u_0\in
D(A)$, if the function $\varphi$ is differentiable on $\mathbb{R}_+$, then LIVP (2.11) has a unique strong solution.
\vskip3mm
The following regularity result will be used in the proof
of our main results.
\vskip3mm
\noindent\textbf{Lemma 2.4} (\cite[Chapter II, Theorem 3.3 ]{te97}). \emph{Assume that $V$ and $H$ are two Hilbert spaces, $V\subset H$, $V$ is dense in $H$, the injection is continuous and compact, $A:D(A)\subset H\to H$ is a positive definite self-adjoint
operator in $H$.  Then for any $u_0\in V$ and $\varphi\in L^2(\mathbb{R}_+,V)$, the mild solution of the
LIVP (2.11) has the regularity
$$u\in  W^{1,2}(\mathbb{R}_+,H)\cap L^2(\mathbb{R}_+,H_1)\cap C(\mathbb{R}_+,V).$$
}\vskip3mm

Denote by
$$C_\omega(\mathbb{R},H)=\{u\mid u:\mathbb{R}\rightarrow H\;\textrm{is continuous and}\;u(t+\omega)=u(t)\; \textrm{for every}\; t\in \mathbb{R}\}.$$
Then it is easy to verify that $C_\omega(\mathbb{R},H)$ is a Banach space endowed with the norm
$$
\| u\|_{C}=\max_{t\in [0,\omega]}\|u(t)\|_2,\quad \forall \;u\in C_\omega(\mathbb{R},H).
$$
\vskip3mm
\noindent\textbf{Lemma 2.5.}\emph{ For every $\varphi\in C_\omega(\mathbb{R},H)$, the linear evolution equation
$$
 u'(t)+Au(t)= \varphi(t),\qquad  t\in \mathbb{R}\eqno(2.13)
$$
has a unique $\omega$-periodic mild solution $u\in C_\omega(\mathbb{R},H)$ which is given by
$$u(t)=\Big(I-T(\omega)\Big)^{-1}\int_{t-\omega}^{t}
T(t-s)\varphi(s)ds,\quad t\in \mathbb{R}.\eqno(2.14)$$
}\vskip2mm
\noindent\textbf{Proof.}\quad By the above discussion,  we know that the evolution equation (2.11) has
a unique mild solution $u$ given by (2.12) and
$$
u(\omega)=T(\omega)u_0+\int_0^{\omega}T(\omega-s)\varphi(s)ds.\eqno(2.15)
$$
From (2.9) one gets that $\|T(\omega)\|\leq e^{-(\gamma\pi^4+\pi^2-1)\omega}<1$. Therefore  we know that $I-T(\omega)$ has a bounded inverse operator $\Big(I-T(\omega)\Big)^{-1}$. Hence, there exists a unique initial value
$$
u_0=\Big(I-T(\omega)\Big)^{-1}\int_0^{\omega}
T(\omega-s)\varphi(s)ds\eqno(2.16)
$$such that the unique mild solution $u$ of LIVP (2.11) expressed by (2.12) satisfies the periodic boundary condition $u(0)=u_0=u(\omega)$. Therefore, from (2.12), (2.15), (2.16) and the fact that  $\varphi(t)=\varphi(t+\omega)$ for $t\in \mathbb{R}$, we get that for every $t\in \mathbb{R}_+$
\begin{eqnarray*}
u(t+\omega)&=&\Big(I-T(\omega)\Big)^{-1}\int_{t}^{t+\omega}
T(t+\omega-s)\varphi(s)ds\\
&=&\Big(I-T(\omega)\Big)^{-1}\int_{t-\omega}^{t}
T(t-s)\varphi(s+\omega)ds=u(t).
\end{eqnarray*}
Therefore, the $\omega$-periodic extension of $u$ on $\mathbb{R}$ is a unique $\omega$-periodic mild solution of linear evolution equation (2.13). Combining  (2.12) and (2.16), we get that the mild solution $u$ of linear evolution equation (2.13) satisfies (2.14).\vskip3mm

Conversely, we can verify directly that the function $u\in C_\omega(\mathbb{R},H)$
given by (2.14) is a mild solution of linear evolution equation (2.13).   This completes the proof of Lemma 2.5. \hfill$\Box$\vskip3mm

In what follows, we recall the Bellman type inequality with delays (see \cite[Lemma 4.1]{l11}), which will be used in the proof of our main results.\vskip3mm
\noindent\textbf{Lemma 2.6.}\emph{ Denote $r=\max\{\tau_1,\tau_2,\cdots,\tau_n\}$. Let $\psi\in C([-r,\infty),\mathbb{R}_+)$. If  there exist positive constants
$b_1$, $b_2$, $\cdots$, $b_n$ such that $\psi$ satisfy the integral inequality
$$
\psi(t)\leq \psi(0)+\sum_{k=1}^{n}b_k\int_0^{t}\psi(s-\tau_k)ds,\quad t\geq 0.
$$
Then for every $t\geq 0$,
$$
\psi(t)\leq \|\psi\|_{C[-r,0]}e^{(\sum_{k=1}^{n}b_k)t},
$$
where $\|\psi\|_{C[-r,0]}=\max\limits_{t\in [-r,0]}|\psi(t)|$.
}
\vskip20mm

\section{Existence and uniqueness of periodic solutions}

\vskip5mm In this section, we will prove the global existence and uniqueness of time $\omega$-periodic classical solutions to the extended Fisher-Kolmogorov equations with delays and general nonlinear term (1.5), i.e., Theorems 1.1 and 1.2.
\vskip5mm\noindent\textbf{Proof of Theorem 1.1.}\quad By the discussions in Section 2, we know that EFK equation (1.5) can be transformed into the abstract delay evolution equation (2.2) in the Hilbert space $H=L^2([0,1],\mathbb{R})$. In what follows, we prove the existence of time $\omega$-periodic mild solutions for abstract delay evolution equation (2.2). Consider
the operator $\mathcal{F}$ on $C_\omega(\mathbb{R},H)$ defined by
\begin{eqnarray*}
\qquad\qquad (\mathcal{F}u)(t)&=&\Big(I-T(\omega)\Big)^{-1}\int_{t-\omega}^{t}
T(t-s)\\[6pt]
& &\cdot[f(u(s-\tau_1),\cdots,u(s-\tau_n))+g(s)]ds,\quad t\in \mathbb{R}.\qquad\qquad\quad(3.1)
\end{eqnarray*}
By the assumptions that $f: \mathbb{R}^n\rightarrow \mathbb{R}$ is locally Lipschitz continuous, $g: [0, 1]\times \mathbb{R}\rightarrow \mathbb{R}$ is H\"{o}lder continuous
and $g(x,t)$ is $\omega$-periodic in $t$ combined with Lemma 2.5 one can easily see that the
operator $\mathcal{F}$ maps $C_\omega(\mathbb{R},H)$ to $C_\omega(\mathbb{R},H)$ is continuous, and the time $\omega$-periodic mild solution of abstract delay evolution equation (2.2) is equivalent to the
fixed point of operator $\mathcal{F}$ defined by (3.1).\vskip2mm

Denote$$
\Omega_R=\{u\in C_\omega(\mathbb{R},H) : \|u\|_C\leq R\},
$$ then $\Omega_R$ is a closed ball in $C_\omega(\mathbb{R},H)$
with center $\theta$ and radius $R$. By the condition (H1) we know  that for any $u\in C_\omega(\mathbb{R},H)$
$$\|f(u(t-\tau_1),\cdots,u(t-\tau_n))+g(t)\|_2\leq\sum_{k=1}^{n}\beta_k\|u(t-\tau_k)\|_2+K,\quad t\in \mathbb{R}.\eqno(3.2)$$
Furthermore, from the fact that $\|T(\omega)\|_2\leq e^{-(\gamma\pi^4+\pi^2-1)\omega}<1$ combined with Neumann series, $\Big(I-T(\omega)\Big)^{-1}$ can be expressed by
$$
\Big(I-T(\omega)\Big)^{-1}=\sum_{n=0}^{\infty}T^n(\omega).
$$ Therefore, by the above equality and (2.9) one gets that
$$
\Big\|\Big(I-T(\omega)\Big)^{-1} \Big\|_2=\Big\|\sum_{n=0}^{\infty}T^n(\omega)\Big\|_2\leq \sum_{n=0}^{\infty}e^{-(\gamma\pi^4+\pi^2-1)n\omega}=
\frac{1}{1-e^{-(\gamma\pi^4+\pi^2-1)\omega}}.\eqno(3.3)
$$
Next, we prove that there exists a constant $R$ big enough such that the operator $\mathcal{F}$ maps $\Omega_R$ to $\Omega_R$. In fact, if we choose
$$
R\geq \frac{K}{\gamma\pi^4+\pi^2-1-\sum_{k=1}^{n}\beta_k},\eqno(3.4)
$$then for any $u\in \Omega_R$ and $t\in \mathbb{R}$, by (2.9), (3.1)-(3.4) and the condition (H2), we
have
\begin{eqnarray*}
\| (\mathcal{F}u)(t)\|_2&\leq&
\frac{1}{1-e^{-(\gamma\pi^4+\pi^2-1)\omega}}
\int_{t-\omega}^{t}e^{-(\gamma\pi^4+\pi^2-1)(t-s)}\Big[\sum_{k=1}^{n}\beta_k\|u(s-\tau_k)\|_2+K\Big]ds
\\[6pt]
&\leq&\frac{1}{1-e^{-(\gamma\pi^4+\pi^2-1)\omega}}\cdot\frac{1-e^{-(\gamma\pi^4+\pi^2-1)\omega}}
{\gamma\pi^4+\pi^2-1}\cdot\Big[\sum_{k=1}^{n}\beta_k\|u\|_C+K\Big]\\[6pt]
&\leq&\frac{1}{\gamma\pi^4+\pi^2-1}\Big(R\sum_{k=1}^{n}\beta_k+K\Big)\\[6pt]
&\leq&R.
\end{eqnarray*}
Therefore,
$$
\|\mathcal{F}u\|_C=\max_{t\in [0,\omega]}\|(\mathcal{F}u)(t)\|_2\leq R,
$$
which means that $\mathcal{F}u\in\Omega_R$. Therefore, we proved that
$\mathcal{F}:\Omega_R\rightarrow\Omega_R$ is a continuous operator.\vskip2mm

Next, we demonstrate that $\mathcal{F}:\Omega_R\rightarrow\Omega_R$ is
a compact operator. To prove this, we
first show that $\{(\mathcal{F}u)(t):u\in \Omega_R\}$ is relatively compact in $H$ for every $t\in \mathbb{R}$. From the periodicity of the operator $(\mathcal{F}u)(t)$ for $t\in \mathbb{R}$ and $u\in \Omega_R$, we only need to prove that $\{(\mathcal{F}u)(t):u\in \Omega_R\}$ is relatively compact in $H$ for  $0\leq t\leq \omega$. From the $\omega$-periodicity of functions $g$ and $u$, it is easy to see that for every $u\in \Omega_R$,
$$
(\mathcal{F}u)(0)=\Big(I-T(\omega)\Big)^{-1}\int_0^{\omega}
T(\omega-s)[f(u(s-\tau_1),\cdots,u(s-\tau_n))+g(s)]ds.\eqno(3.5)
$$
For any $0<\epsilon<\omega$ and $u\in \Omega_R$, we define the operator $\mathcal{F}^{\,\epsilon}_0$ by
\begin{eqnarray*}
(\mathcal{F}^{\,\epsilon}_0u)(0)&=&\Big(I-T(\omega)\Big)^{-1}\int_0^{\omega-\epsilon}
T(\omega-s)[f(u(s-\tau_1),\cdots,u(s-\tau_n))+g(s)]ds\\[6pt]
&=&T(\epsilon)\Big(I-T(\omega)\Big)^{-1}\int_0^{\omega-\epsilon}
T(\omega-s-\epsilon)\\[6pt]
& & \cdot[f(u(s-\tau_1),\cdots,u(s-\tau_n))+g(s)]ds.\qquad\qquad\qquad\qquad\qquad\qquad (3.6)
\end{eqnarray*}
Since $T(t)$ is compact for every $t>0$, the set $\{(\mathcal {F}^{\,\epsilon}_0u)(0):u\in \Omega_R\}$ is relatively compact in $H$ for every $\epsilon\in (0,\omega)$. Moreover, for every $u\in \Omega_R$, by (3.2), (3.3), (3.5) and (3.6), we get that
\begin{eqnarray*}
 & &\|(\mathcal {F} u)(0)-(\mathcal {F}^{\,\epsilon}_0u)(0)\|_2\\[6pt]
&=&\Big\|\Big(I-T(\omega)\Big)^{-1}\int_{\omega-\epsilon}^{\omega}
T(\omega-s)[f(u(s-\tau_1),\cdots,u(s-\tau_n))+g(s)]ds\Big\|_2\\[6pt]
 &\leq& \frac{1}{1-e^{-(\gamma\pi^4+\pi^2-1)\omega}}\int_{\omega-\epsilon}^{\omega}
 e^{-(\gamma\pi^4+\pi^2-1)(\omega -s)}\Big[\sum_{k=1}^{n}\beta_k\|u(s-\tau_k)\|_2+K\Big]ds\\[6pt]
 &\leq&\frac{1-e^{-(\gamma\pi^4+\pi^2-1)\epsilon}}{(\gamma\pi^4+\pi^2-1)(1-e^{-(\gamma\pi^4+\pi^2-1)\omega})}\Big(R\sum_{k=1}^{n}\beta_k+K\Big)\\[6pt]
&\rightarrow&0\quad \textrm{as}\quad \epsilon\rightarrow 0.
\end{eqnarray*}
Therefore, we have proved that there exists a relatively compact set $\{(\mathcal {F}^{\,\epsilon}_0u)(0):u\in \Omega_R\}$ arbitrarily close to the set $\{(\mathcal {F}u)(0):u\in \Omega_R\}$, this means that the set $\{(\mathcal {F} u)(0):u\in \Omega_R\}$ is relatively compact in $H$. Let $0<t\leq \omega$ be given, $0<\epsilon<t$ and $u\in \Omega_R$, we define the operator $\mathcal {F}^{\,\epsilon}u$ by
\begin{eqnarray*}
 (\mathcal {F}^{\,\epsilon}u)(t)&=&T(t)(\mathcal {F} u)(0)+\int_0^{t-\epsilon}T(t-s)[f(u(s-\tau_1),\cdots,u(s-\tau_n))+g(s)]ds\\[6pt]
&=&T(t)(\mathcal {F} u)(0)\\[6pt]
& &+T(\epsilon)\int_0^{t-\epsilon}T(t-s-\epsilon)[f(u(s-\tau_1),\cdots,u(s-\tau_n))+g(s)]ds.\qquad (3.7)
\end{eqnarray*}
By compactness of the operator $ T(t)$ for $t>0$ combined with the fact that the set $\{(\mathcal {F} u)(0):u\in \Omega_R\}$ is relatively compact in $H$, the set $\{(\mathcal {F}^{\,\epsilon}u)(t):u\in \Omega_R\}$ is relatively compact in $H$ for every $\epsilon\in (0,t)$ and $0<t\leq \omega$. Furthermore, for every $u\in \Omega_R$, by (3.1), (3.2) and (3.7), we get that

\begin{eqnarray*}
\|(\mathcal {F} u)(t)-(\mathcal {F}^{\,\epsilon}u)(t)\|_2&=&\Big\|\int_{t-\epsilon}^{t}
T(t-s)[f(u(s-\tau_1),\cdots,u(s-\tau_n))+g(s)]ds\Big\|_2\\[6pt]
 &\leq& \int_{t-\epsilon}^{t}
 e^{-(\gamma\pi^4+\pi^2-1)(t -s)}\Big[\sum_{k=1}^{n}\beta_k\|u(s-\tau_k)\|_2+K\Big]ds\\[6pt]
 &\leq&\frac{1-e^{-(\gamma\pi^4+\pi^2-1)\epsilon}}{\gamma\pi^4+\pi^2-1}\Big(R\sum_{k=1}^{n}\beta_k+K\Big)\\[6pt]
&\rightarrow&0\quad \textrm{as}\quad \epsilon\rightarrow 0.
\end{eqnarray*}
Hence, we have proved that there exists relatively compact set
$\{(\mathcal {F}^{\,\epsilon}u)(t):u\in \Omega_R\}$ arbitrarily close to the set $\{(\mathcal {F} u)(t):u\in \Omega_R\}$ in $H$ for $0<t\leq\omega$. Therefore,
the set $\{(\mathcal {F} u)(t):u\in \Omega_R\}$ is also relatively compact in $H$ for $0<t\leq\omega$, which combined with the fact that the set $\{(\mathcal {F} u)(0):u\in \Omega_R\}$ is relatively compact in $H$ we get the relatively compactness of the set $\{(\mathcal {F} u)(t):u\in \Omega_R\}$ in $H$ for $0\leq t\leq\omega$.\vskip2mm

In the following we prove that $\{\mathcal{F}(u):u\in\Omega_R\}$ is an equicontinuous family. For any $u\in \Omega_R$ and $t_1$, $t_2\in \mathbb{R}$ with $t_1<t_2$, we get form (2.9), (3.1) and (3.2) that
\begin{eqnarray*}
\|(\mathcal{F}u)(t_2)-(\mathcal{F}u)(t_1)\|_2&\leq& \Big\|(T(t_2)-T(t_1))\Big(I-T(\omega)\Big)^{-1}\int_0^{\omega}
T(\omega-s)\\[6pt]&  & \cdot[f(u(s-\tau_1),\cdots,u(s-\tau_n))+g(s)]ds\Big\|_2\\[6pt]
& &+\Big(R\sum_{k=1}^{n}\beta_k+K\Big)\int_{t_1}^{t_2}
e^{-(\gamma\pi^4+\pi^2-1)(t_2-s)}ds\qquad\qquad\qquad
\end{eqnarray*}
\begin{eqnarray*}\qquad\qquad\qquad\qquad& &+\Big(R\sum_{k=1}^{n}\beta_k+K\Big)\int_{0}^{t_1}
\|T(t_2-s)-T(t_1-s)\|ds\\[6pt]
&:=&I_1+I_2+I_3,
\end{eqnarray*}
where
$$I_1=\Big\|(T(t_2)-T(t_1))\Big(I-T(\omega)\Big)^{-1}\int_0^{\omega}
T(\omega-s)[f(u(s-\tau_1),\cdots,u(s-\tau_n))+g(s)]ds\Big\|_2,$$
$$I_2=\Big(R\sum_{k=1}^{n}\beta_k+K\Big)\int_{t_1}^{t_2}
e^{-(\gamma\pi^4+\pi^2-1)(t_2-s)}ds,\qquad\qquad\qquad\qquad\qquad\qquad\qquad\qquad\qquad$$
$$I_3=\Big(R\sum_{k=1}^{n}\beta_k+K\Big)\int_{0}^{t_1}
\|T(t_2-s)-T(t_1-s)\|ds.\qquad\qquad\qquad\qquad\qquad\qquad\qquad\qquad\qquad$$
Therefore, we only need to check $I_i$ tend to
$0$ independently of $u\in \Omega_R$ when $t_2-t_1\rightarrow 0$ for
$i=1,2,3$. For $I_1$, by the definition of $I_1$, (2.9), (3.2) and (3.3), we get that
\begin{eqnarray*}& &\Big\|\Big(I-T(\omega)\Big)^{-1}\int_0^{\omega}
T(\omega-s)[f(u(s-\tau_1),\cdots,u(s-\tau_n))+g(s)]ds\Big\|_2\\[6pt]
 &\leq&\frac{1}{1-e^{-(\gamma\pi^4+\pi^2-1)\omega}}
\int_0^{\omega}e^{-(\gamma\pi^4+\pi^2-1)(\omega-s)}\Big[\sum_{k=1}^{n}\beta_k\|u(s-\tau_k)\|_2+K\Big]ds\\
 &\leq& \frac{R\sum\limits_{k=1}^{n}\beta_k+K}{\gamma\pi^4+\pi^2-1}.
\end{eqnarray*}
From the above inequality combined with the strong continuity of the
semigroup $T(t)$ $(t\geq0)$ and the definition of $I_1$, we can easily get that $I_1\rightarrow0$ as $t_2-t_1\rightarrow 0$. For $I_2$, we can get by direct calculus that
\begin{eqnarray*}I_2&\leq&
\frac{R\sum\limits_{k=1}^{n}\beta_k+K}{\gamma\pi^4+\pi^2-1}
\Big[1-e^{-(\gamma\pi^4+\pi^2-1)(t_2-t_1)}\Big]\\[6pt]
&\rightarrow&0\quad \textrm{as}\quad
t_2-t_1\rightarrow0.\end{eqnarray*}
For $I_3$, by the definition of $I_3$, the Lebesgue dominated convergence theorem as well as the norm continuity and uniform boundedness of $T(t)$ for $t>0$, we get that
\begin{eqnarray*}
I_3&\leq&
\Big(R\sum_{k=1}^{n}\beta_k+K\Big)\int_{0}^{t_1}
\|T(t_2-t_1+s)-T(s)\|ds
\\[6pt]
&\rightarrow&0\quad \textrm{as}\quad t_2-t_1\rightarrow 0.
\end{eqnarray*}
As a result, $\|(\mathcal {F} u)(t_2)-(\mathcal {F} u)(t_1)\|_2$ tends to zero
independently of $u\in \Omega_R$ as $t_2-t_1\rightarrow 0$, which
means that the family $\{\mathcal{F}(u):u\in\Omega_R\}$ is equicontinuous. Therefore, $\{\mathcal{F}u:u\in \Omega_R\}$ is relatively compact by Arzela-Ascoli theorem. Hence, the continuity of the operator $\mathcal{F}$ and the relatively compactness of the set $\{\mathcal{F}u:u\in \Omega_R\}$ imply that  $\mathcal{F}:\Omega_R \rightarrow\Omega_R$ is a completely continuous operator. It follows from
Schauder's fixed point theorem that $\mathcal{F}$ has at least one fixed point $u\in \Omega_R$, which is just a time $\omega$-periodic mild solution of the abstract delay evolution equation (2.2). \vskip 2mm

In what follows, we prove the regularity for the time $\omega$-periodic mild solution $u$ of the abstract delay evolution equation (2.2).
Since $u$ is the mild solution of the linear evolution equation (2.13)
for $\varphi(\cdot)=f(u(\cdot-\tau_1),\cdots,u(\cdot-\tau_n))+g(\cdot)\in L^2(\mathbb{R},H)$,
by the maximal regularity of linear evolution equations with
positive definite operator in Hilbert spaces (see for details Lemma 2.4), when $u_0\in V:=H_{\frac{1}{2}}$, the mild solution of
LIVP (2.11) has the regularity
$$u\in  W^{1,2}(\mathbb{R}_+,H)\cap L^2(\mathbb{R}_+,H_1)\cap C(\mathbb{R}_+,H_{\frac{1}{2}})\eqno(3.8)$$
and it is a strong solution. We notice that $u(t)$ is the mild solution of
LIVP (2.11) for
$$u_0=\Big(I-T(\omega)\Big)^{-1}\int_0^{\omega}
T(\omega-s)\varphi(s)ds.$$ By the
representation (2.12) of mild solution, $u(t)=T(t)u_0+v(t)$, where
$v(t)=\int_0^{t}T(t-s)\varphi(s)ds$. Since the function $v(t)$ is a mild
solution of LIVP (2.11) with the null initial value
$u(0)=\theta$, $v$ has the regularity (3.8). By the analytic
property of the semigroup $T(t)$, $T(\omega)u_0\in D(A)\subset
H_{1/2}$. Hence,
$$u_0=u(\omega)=T(\omega)u_0+v(\omega)\in
H_{1/2}.$$ Using the regularity (3.8) again, we obtain that $u\in
W^{1,2}(\mathbb{R},H)\cap L^2(\mathbb{R},D(A))$ and it is a  time $\omega$-periodic strong solution of the linear evolution equation (2.13), which means that the fixed point $u$ of the operator $\mathcal{F}$ defined by (3.1) belongs to $W^{1,2}(\mathbb{R},H)\cap L^2(\mathbb{R},D(A))$ is the time $\omega$-periodic strong solution of the abstract delay evolution equation (2.2). Furthermore, by the usual regularization
method via the theory of analytic semigroups of linear operators used in \cite[Lemma 4.2]{a78} combined with the fact that $f: \mathbb{R}^n\rightarrow \mathbb{R}$ is locally Lipschitz continuous and $g: \mathbb{R}\times[0, 1]\rightarrow \mathbb{R}$ is H\"{o}lder continuous, we can prove that $u\in C^{1,2}(\mathbb{R}\times[0,1])$ is a time $\omega$-periodic classical solution of EFK equation (1.5). This completes the proof of Theorem 1.1.\hfill$\Box$\vskip5mm

\noindent\textbf{Proof of Theorem 1.2.}\quad
By the proof of Theorem 1.1 we know that EFK equation (1.5) can be transformed into the abstract delay evolution equation (2.2) in Hilbert space $H=L^2([0,1],\mathbb{R})$ and the time $\omega$-periodic mild solutions of abstract delay evolution equation (2.2) is equivalent to the
fixed point of operator $\mathcal{F}$ defined by (3.1), which  maps $C_\omega(\mathbb{R},H)$ to $C_\omega(\mathbb{R},H)$. By the condition (H3) we know that for every $t\in \mathbb{R}$ and $u$, $v\in C_\omega(\mathbb{R},H)$
\begin{eqnarray*}
\qquad\|f(u(t-\tau_1),\cdots,u(t-\tau_n))&-&f(v(t-\tau_1),\cdots,v(t-\tau_n))\|_2\\[6pt]
&\leq&\sum_{k=1}^{n}\beta_k\|u(t-\tau_k)-v(t-\tau_k)\|_2.\qquad\qquad\qquad(3.9)
\end{eqnarray*}
Therefore, for any $u$, $v\in C_\omega(\mathbb{R},H)$, by (2.9), (3.1), (3.3) and (3.9), we get that
\begin{eqnarray*}
& &\| (\mathcal{F}u)(t)- (\mathcal{F}v)(t)\|_2\\[6pt]&\leq&
\frac{1}{1-e^{-(\gamma\pi^4+\pi^2-1)\omega}}
\int_{t-\omega}^{t}e^{-(\gamma\pi^4+\pi^2-1)(t-s)}\Big(\sum_{k=1}^{n}\beta_k\|u(s-\tau_k)-v(s-\tau_k)\|_2\Big)ds
\\[6pt]
&\leq&\frac{1}{1-e^{-(\gamma\pi^4+\pi^2-1)\omega}}\cdot\frac{1-e^{-(\gamma\pi^4+\pi^2-1)
\omega}}{\gamma\pi^4+\pi^2-1}\cdot\Big(\sum_{k=1}^{n}\beta_k\|u-v\|_C\Big)\\[6pt]
&=&\frac{\sum_{k=1}^{n}\beta_k}{\gamma\pi^4+\pi^2-1}\|u-v\|_C,\end{eqnarray*}
which means that,
$$
\|\mathcal{F}u-\mathcal{F}v\|_C=\max_{t\in [0,\omega]}\|(\mathcal{F}u)(t)-(\mathcal{F}v)(t)\|_2\leq \frac{\sum_{k=1}^{n}\beta_k}{\gamma\pi^4+\pi^2-1}\|u-v\|_C.
$$
From the condition (H2) it follows that $\mathcal{F}:C_\omega(\mathbb{R},H)\rightarrow C_\omega(\mathbb{R},H)$ is a contraction operator, and therefore $\mathcal{F}$ has a
unique fixed point $u\in C_\omega(\mathbb{R},H)$, which is in turn the unique time $\omega$-periodic mild solution of the abstract delay evolution equation (2.2). By using a completely similar method with which used in the proof of Theorem 1.1 combined with the fact that $f: \mathbb{R}^n\rightarrow \mathbb{R}$ is locally Lipschitz continuous and $g: \mathbb{R}\times[0, 1]\rightarrow \mathbb{R}$ is H\"{o}lder continuous, we can prove that $u\in C^{1,2}(\mathbb{R}\times[0,1])$ is the unique time $\omega$-periodic classical solution of EFK equation (1.5). This completes the proof of Theorem
1.2.\hfill$\Box$\vskip20mm

\section{Global asymptotic stability of periodic solutions}
\vskip3mm In this section, we will prove the global asymptotic stability of time $\omega$-periodic classical solution for EFK equation (1.5), i.e., Theorem 1.3.
For this purpose, we firstly discuss the existence of classical solutions to the initial value problem for extended Fisher-Kolmogorov equations with delays
$$
\left\{\begin{array}{ll}
\frac{\partial}{\partial t}u(t,x)+\gamma \frac{\partial^4}{\partial x^4}u(t,x)-\frac{\partial^2 }{\partial x^2}u(t,x)-u(t,x)=g(t,x)\\[6pt]
\qquad\qquad\qquad\qquad\quad+f(u(t-\tau_1,x),\cdots,u(t-\tau_n,x)),\;\textrm{in}\; \mathbb{R}_+\times (0,1),\\[12pt]
u(t,0)=u(t,1)=0,\quad u_{xx}(t,0)=u_{xx}(t,1)=0,\quad t\in \mathbb{R}_+,\\[12pt]
u(t,x)=\kappa(t,x),\qquad t\in [-r,0], \quad x\in (0,1),
 \end{array} \right.\eqno(4.1)
$$
where  $\gamma>0$ is a constant, $f: \mathbb{R}^n\rightarrow \mathbb{R}$ is a nonlinear continuous function, $g: \mathbb{R}_+\times[0, 1]\rightarrow \mathbb{R}$ is continuous, $\tau_1$, $\tau_2$, $\cdots$, $\tau_n$ are positive constants, $r=\max\{\tau_1,\tau_2,\cdots,\tau_n\}$, $\kappa\in C([-r,0]\times(0,1), \mathbb{R})$.

Let $\kappa(t)=\kappa(t,\cdot)$ for $t\in [-r,0]$. Then from the discussion in Section 2 we know that the initial value problem of extended Fisher-Kolmogorov  equations with delays  (4.1) can be transformed into the abstract form of initial value problem to delay evolution equation
$$
\left\{\begin{array}{ll}
  u'(t)+ Au(t)=f(u(t-\tau_1),\cdots,u(t-\tau_n))+g(t),\quad t\in \mathbb{R}_+, \\[12pt]
   u(t)=\kappa(t)\qquad t\in [-r,0]
    \end{array} \right.\eqno(4.2)
$$
in the Hilbert space $H=L^2((0,1),\mathbb{R})$. A function $u\in C([-r,\infty),H)$ is said to be a mild solution of initial value problem (4.2) if $u(t)$ satisfies
$$
u(t)=T(t)u(0)+\int_0^{t}T(t-s)[f(u(s-\tau_1),\cdots,u(s-\tau_n))+g(s)]ds\quad\textrm{for}\quad t\in\mathbb{R}_+, \eqno(4.3)
$$
and the initial condition
$$u(t)=\kappa(t)\quad\textrm{for}\quad t\in [-r,0].\eqno(4.4)$$
\vskip5mm\noindent\textbf{Theorem 4.1.} \emph{Assume that $f: \mathbb{R}^n\rightarrow \mathbb{R}$ is locally Lipschitz continuous and $g: \mathbb{R}_+\times[0, 1]\rightarrow \mathbb{R}$ is H\"{o}lder continuous. If the conditions (H2) and (H3) are satisfied, then the initial value problem of extended Fisher-Kolmogorov  equations with delays  (4.1) has a unique  classical solution $u\in C^{1,2}(\mathbb{R}\times[0,1])$.}
\vskip3mm\noindent\textbf{Proof.}\quad By the above discussion, we know that the initial value problem of extended Fisher-Kolmogorov  equations with delays  (4.1) can be transformed into the abstract form of initial value problem to delay evolution equation (4.2) in Hilbert space $H=L^2([0,1],\mathbb{R})$.  Define
the operator $\mathcal{Q}$ on $ C([-r,\infty),H)$ as follows
$$
(\mathcal{Q}u)(t)=\left\{\begin{array}{ll}
  T(t)u(0)+\int_0^{t}
T(t-s)[f(u(s-\tau_1),\cdots,u(s-\tau_n))+g(s)]ds,\; t\in \mathbb{R}_+,\\[12pt]
  \kappa(t),\qquad t\in [-r,0]
    \end{array} \right.\eqno(4.5)
$$
Then by the assumptions that $f: \mathbb{R}^n\rightarrow \mathbb{R}$ is locally Lipschitz continuous, $g: \mathbb{R}_+\times[0, 1]\rightarrow \mathbb{R}$ is H\"{o}lder continuous
and  $\kappa \in C([-r,0]\times(0,1), \mathbb{R})$ one can easily see that $\mathcal{Q}$ maps $C([-r,\infty),H)$ to $ C([-r,\infty),H)$ and the mild solutions of initial value problem for delay evolution equation (4.2) is equivalent to the
fixed point of operator $\mathcal{Q}$ defined by (4.5).

For any $u$, $v\in C([-r,\infty),H)$, by (2.9), (3.9) and (4.5), we get that
\begin{eqnarray*}
\| (\mathcal{Q}u)(t)- (\mathcal{Q}v)(t)\|_2&\leq&
\int_0^{t}
e^{-(\gamma\pi^4+\pi^2-1)(t-s)}\Big(\sum_{k=1}^{n}\beta_k\|u(s-\tau_k)-v(s-\tau_k)\|_2\Big)ds\\[6pt]
&\leq&\frac{\sum_{k=1}^{n}\beta_k}{\gamma\pi^4+\pi^2-1}\|u-v\|_C,
\end{eqnarray*}
which means that,
$$
\|\mathcal{Q}u-\mathcal{Q}v\|_C=\sup_{t\in [-r,\infty)}\|(\mathcal{Q}u)(t)-(\mathcal{Q}v)(t)\|_2\leq\frac{\sum_{k=1}^{n}\beta_k}{\gamma\pi^4+\pi^2-1} \|u-v\|_C.
$$
From the condition (H2) it follows that $\mathcal{Q}:C([-r,\infty),H)\rightarrow C([-r,\infty),H)$ is a contraction operator, and therefore $\mathcal{Q}$ has a
unique fixed point $u\in C([-r,\infty),H)$, which is in turn the unique  mild solution of the initial value problem to delay evolution equation (4.2). By using a completely similar method to the one used in the proof of Theorem 1.1 combined with the fact that $f: \mathbb{R}^n\rightarrow \mathbb{R}$ is locally Lipschitz continuous and $g: \mathbb{R}\times[0, 1]\rightarrow \mathbb{R}$ is H\"{o}lder continuous, we can prove that $u\in C^{1,2}(\mathbb{R}_+\times[0,1])$ is the unique classical solution  for the initial value problem of extended Fisher-Kolmogorov  equations with delays  (4.1). This completes the proof of Theorem 4.1.\hfill$\Box$\vskip5mm

\noindent\textbf{Proof of Theorem 1.3.}\quad One can easily see that (H2)$'$ $\Rightarrow$ (H2). Therefore, By Theorem 1.2 we know that EFK equation (1.5) has a unique time $\omega$-periodic classical solution $\overline{u}\in C^{1,2}(\mathbb{R}\times[0,1])$. Furthermore, Theorem 4.1 means that the initial value problem of extended Fisher-Kolmogorov  equations with delays  (4.1) has a unique classical solution $u_\kappa\in C^{1,2}(\mathbb{R}_+\times[0,1])$.

By (2.9), (3.1), (3.5), (3.9) and (4.5), we get that
\begin{eqnarray*}
& &\|\overline{u}(t)- u_\kappa(t)\|_2\\[6pt]&\leq&
e^{-(\gamma\pi^4+\pi^2-1)t}\|\overline{u}(0)- u_\kappa(0)\|_2
\\[6pt]& &+\int_0^{t}
e^{-(\gamma\pi^4+\pi^2-1)(t-s)}\Big(\sum_{k=1}^{n}\beta_k\|\overline{u}(s-\tau_k)-u_\kappa(s-\tau_k)\|_2\Big)ds\\[6pt]
&=&e^{-(\gamma\pi^4+\pi^2-1)t}\|\overline{u}(0)- u_\kappa(0)\|_2\\[6pt]
& &+e^{-(\gamma\pi^4+\pi^2-1)t}\sum_{k=1}^{n}\beta_ke^{(\gamma\pi^4+\pi^2-1)\tau_k}\int_0^{t}
e^{(\gamma\pi^4+\pi^2-1)(s-\tau_k)}\|\overline{u}(s-\tau_k)-u_\kappa(s-\tau_k)\|_2ds,
\end{eqnarray*}
from which one gets that
\begin{eqnarray*}
& &e^{(\gamma\pi^4+\pi^2-1)t}\|\overline{u}(t)- u_\kappa(t)\|_2\\[6pt]&\leq&
\|\overline{u}(0)- u_\kappa(0)\|_2\\[6pt]
& &+\sum_{k=1}^{n}\beta_ke^{(\gamma\pi^4+\pi^2-1)\tau_k}\int_0^{t}
e^{(\gamma\pi^4+\pi^2-1)(s-\tau_k)}\|\overline{u}(s-\tau_k)-u_\kappa(s-\tau_k)\|_2ds\qquad\qquad (4.6)
\end{eqnarray*}
Letting
$$
\psi(t)=e^{(\gamma\pi^4+\pi^2-1)t}\|\overline{u}(t)- u_\kappa(t)\|_2,\quad [-r,\infty).
$$
Then from (4.6) we get that
$$
\psi(t)\leq \psi(0)+\sum_{k=1}^{n}\beta_ke^{(\gamma\pi^4+\pi^2-1)\tau_k}\int_0^{t}
\psi(s-\tau_k)ds,\qquad t\geq0.\eqno(4.7)
$$
Therefore, by (4.7) and Lemma 2.6, we  know that for every $t\geq 0$,
\begin{eqnarray*}
e^{(\gamma\pi^4+\pi^2-1)t}\|\overline{u}(t)- u_\kappa(t)\|_2&=&\psi(t)\leq \max\limits_{t\in [-r,0]}e^{(\gamma\pi^4+\pi^2-1)t}\|\overline{u}(t)- \kappa(t)\|_2\\[6pt]
& &\cdot e^{[\sum_{k=1}^{n}\beta_k e^{(\gamma\pi^4+\pi^2-1)\tau_k}]t},
\end{eqnarray*}
from which one gets that
$$
\|\overline{u}(t)- u_\kappa(t)\|_2\leq \max\limits_{t\in [-r,0]}e^{(\gamma\pi^4+\pi^2-1)t}\|\overline{u}(t)- \kappa(t)\|_2e^{[\sum_{k=1}^{n}\beta_k e^{(\gamma\pi^4+\pi^2-1)\tau_k}-(\gamma\pi^4+\pi^2-1)]t}.\eqno(4.8)
$$
By the condition (H2)$'$, we get that $$\sum\limits_{k=1}^{n}\beta_k e^{(\gamma\pi^4+\pi^2-1)\tau_k}-(\gamma\pi^4+\pi^2-1)<0.\eqno(4.9)$$
Hence, from (4.8) and (4.9) we know that
$$
\Big(\int_0^1|\overline{u}(t,x)- u_\kappa(t,x)|^2dx\Big)^{\frac{1}{2}}=\|u(t)- u_\kappa(t)\|_2\rightarrow 0\quad \textrm{as} \quad t\rightarrow +\infty.
$$
Therefore, the $\omega$-periodic classical solution $u$ of EFK equation (1.5) is globally asymptotically stable and it exponentially
attracts every classical solution for the initial value problem of extended Fisher-Kolmogorov  equations with delays.
This completes the proof of Theorem 1.3.\hfill$\Box$

\section*{Acknowledgments} The authors would like to express sincere thanks to the anonymous referee for his/her carefully reading the manuscript and valuable comments
and suggestions.

\end{document}